\documentclass[12pt]{article}

\usepackage{amsmath, amssymb, latexsym, amscd, amsthm}

\newtheorem{theorem}{Theorem}[section]

\newtheorem{lemma}[theorem]{Lemma}

\newtheorem{statement}[theorem]{Statement}

\theoremstyle{definition}
\newtheorem{definition}[theorem]{Definition}
\newtheorem{remark}[theorem]{Remark}

\def\diam{\operatorname{diam}}
\def\id{\operatorname{id}}

\begin{document}


\title{Wild high-dimensional
Cantor fences in $\mathbb R^n$, Part I}   
\author{Olga Frolkina\footnote{Supported by Russian Foundation of Basic Research; Grant
No.~15--01--06302.}\\
Chair of General Topology and Geometry,\\
Faculty of Mechanics and Mathematics,\\
M.V.~Lomonosov Moscow State University,\\
Leninskie Gory 1, GSP-1,\\
Moscow 119991, Russia\\
E-mail: olga-frolkina@yandex.ru
}

\maketitle

\begin{abstract}
Let $\mathcal C$ be the Cantor set.
For each $n\geqslant 3$ 
we construct an embedding
$A: \mathcal C \times \mathcal C \to \mathbb R^n$
such that
$A(\mathcal C \times \{s\})$, for $s\in\mathcal C$,
are
pairwise ambiently incomparable
everywhere wild Cantor sets
(generalized Antoine's necklaces).
This serves as a base for another new result 
proved in this paper:
for each $n\geqslant 3$ 
and
any non-empty perfect compact set
$X$ which is embeddable in $\mathbb R ^{n-1}$,
we describe
an embedding
$\mathbb A : X \times \mathcal C \to \mathbb R^n$
such that
each
$\mathbb A (X \times \mathcal \{s\} )$,
$s\in \mathcal C$,
contains the corresponding
$A (\mathcal C \times \{s\} )$,
and is ``nice''
on the complement 
$\mathbb A (X \times \mathcal \{s\} )-A (\mathcal C \times \{s\} )$;
in particular,
the images
$\mathbb A ( X \times \{s\})$,
for $s\in\mathcal C$, 
are 
ambiently incomparable pairwise 
disjoint copies of $X$.
This generalizes and strengthens
theorems of J.R.~Stallings (1960), R.B.~Sher (1968), and B.L.~Brechner--J.C.~Mayer
(1988). 


\end{abstract}

\section{Introduction and
Statements of Main Results}

Theory of wild embeddings comes back
to papers of L.~Antoine, P.S.~Urysohn, J.W.~Alexander.
Let us refer to 
the reviews
\cite{BC},
\cite{Burgess-75}, \cite{Daverman}, \cite{Cernavsky}
and the
books \cite{Rushing},
\cite{DV}
which
contain hundreds of references on the subject.

It is known that for $n\geqslant 3$
there exist uncountably  many
inequivalent embeddings of the Cantor set
in $\mathbb R^n$
\cite{Sher68}, \cite[Cor.~5.2]{Wright-rigid};
this together with
``feelers method''
(see Statement~\ref{feelers-ext})
implies that 
each non-empty perfect compact
subset of $S^{n-1}$
can be embedded in $\mathbb R^n$ in uncountably many
inequivalent ways, compare
\cite[p.~114]{Sher69},
\cite[Thm.~7.1]{Wright-rigid}, 
\cite[Thm.~2, 3, and Remark on p.~320]{BM}.

In this paper, we study {\it disjoint\/} wild embeddings.
We would like to
embed uncountably many
copies of a compactum simultaneously,
so that they are mutually exclusive.
The assumption of being wild is essential.
Concentric spheres
of arbitrary radii
form a family of cardinality continuum; 
in contrast to this,
by results of R.H.~Bing, 
it is impossible to place
an uncountable collection
of pairwise disjoint
wild closed surfaces in $\mathbb R^3$ 
(scheme of proof
is given in \cite{Bing-abstr};
the full proof needs results of
\cite{Bing59} and \cite{Bing-TAMS61};
see also 
a short sketch of Bing's idea in
\cite[Thm.~3.6.1]{BC}).
For higher dimensions,
Bing's 
non-embedding result 
has 
at the moment
only partial
generalizations, see
\cite[Thm. 1, 2]{Bryant},
\cite[Thm. 10.5]{Burgess-75},
\cite[p.383, Thm. 3C.2]{Daverman}.

J.R.~Stallings
constructed a
family
of continuum cardinality of
pairwise disjoint wild 2-disks in
$\mathbb R^3$ \cite{Stallings}.
J.~Martin showed
that
all except countably many disks in such a
collection must be locally tame except on their 
boundaries, hence must lie on 2-spheres \cite{Martin}
(this result was predicted in
\cite{Bing-TAMS61} where important facts needed
for Martin's arguments are proved).

R.B.~Sher 
modified Stallings' construction so that 
no two
disks of the family 
are ambiently homeomorphic, that
is, no self-homeomorphism of $\mathbb R^3$ 
can map one disk onto another
\cite{Sher-note}.

J.L.~Bryant
noticed that the
Stallings' construction can be extended
to the case of $(n-1)$-disks in $\mathbb R^n$, $n\geqslant 4$
by taking direct product with $I^{n-3}$
\cite[p.~479]{Bryant}.
Our construction is more
complicated, and it provides 
disjoint embeddings
with additional property
of pairwise incomparability,
not only for disks, but
for arbitrary perfect compacta.
Our method unifies and extends ideas of
\cite{Stallings}, \cite{Sher-note},
\cite{BM}.

\begin{theorem}\label{Thm2}
For each $n\geqslant 3$ 
and
any non-empty perfect compact set
$X$ which is embeddable in $\mathbb R^{n-1}$,
there exists
an embedding
$\mathbb A : X \times \mathcal C \to \mathbb R^n$
such that
the images
$\mathbb A ( X \times \{s\})$,
for $s\in\mathcal C$, 
are
non-locally-hyperplanar and
pairwise
ambiently incomparable 
copies of $X$ in $\mathbb R^n$.

Further:

(i) $\mathbb A$ extends to an isotopy
$X\times I \to \mathbb R^n$;

(ii) if, in addition, $X$ is a polyhedron,
then each $\mathbb A(X\times \{s\})$,
$s\in \mathcal C$, is wild.
\end{theorem}

The notion of local hyperplanarity is
introduced in Definition~\ref{hyperplanar},
see below.

\begin{remark}
Let us make some remarks concerning possible analogues and generalizations of Theorem \ref{Thm2}.

1) There is no evident analogue
for countable, compact sets:
any countable, compact set
can be straightened by an ambient
isotopy
\cite[1.1, 1.2]{Ivanov-diss},
see also
\cite[Thm.~I.4.2]{Keldysh}.

2)
By 
\cite[Thm.~1.8.10, 4.1.5]{Engelking}, each $n$-dimensional
compact set in $\mathbb R^n$
contains an open
ball; hence
Theorem \ref{Thm2} does not 
extend to the case of
$n$-dimensional
compacta $X$.

3) 
There is no evident analogue of Theorem \ref{Thm2}
for $n=2$.
In $\mathbb  R^2$ each zero-dimensional compactum
is tame 
\cite[{\bf 75}, p.~87--89]{Antoine-diss},
\cite[Cor.~II.3.2, Cor.~II.3.3]{Keldysh},
\cite[Chap.~13]{Moise}.
E.D.~Tymchatyn and R.B.~Walker showed 
that for
each embedding
$F:I\times \mathcal C \hookrightarrow
\mathbb R^2$ 
there exists a
homeomorphism $h:\mathbb R^2\cong\mathbb R^2$
such that
$h\circ F(I\times \mathcal C)=I\times \mathcal C$~\cite{TW}.
\end{remark}

Our proof of Theorem \ref{Thm2} is based on the following two theorems which are interesting in themselves:

\begin{theorem}\label{MainLemma1}
Let $n\geqslant 3$.
There exists
an embedding
$A: \mathcal C \times \mathcal C \to \mathbb R^n$
such that
for each $s\in \mathcal C$
the corresponding Cantor set
$A( \mathcal C  \times \{ s\} ) = : \mathcal A_s$
is everywhere wild;
moreover, for $s\neq t$, the 
sets
$\mathcal A_s$ and $\mathcal A_t$ are
ambiently incomparable.
In fact,  
each $\mathcal A_s$
is an Antoine's necklace for $n=3$ and a 
generalized Antoine's necklace for $n\geqslant 4$.
\end{theorem}

\begin{theorem}\label{reembedding}
For 
any non-empty perfect compact set
$X\subset \mathbb R^{n}$, where $n\geqslant 2$,
there exists
a homeomorphism
$F: \mathbb R^n \to \mathbb R^{n}$
such that
$F(X)\subset I^n$ and
$F(X) \cap \partial I^n
= \{ 0\}^{n-1}\times\mathcal C
$,
where
$\mathcal C\subset I$ is the middle-thirds Cantor set.
\end{theorem}

Before we finish the introductory part,
let us mention some other results in this field.

R.B.~Sher (assuming the Continuum Hypothesis)
constructed 
an uncountable
collection of pairwise disjoint arcs in $\mathbb R^3$
which contain a representative
(under self-homeomorphism of $\mathbb R^3$)
of each arc locally tame modulo
a compact 0-dimensional set \cite[Thm.~3]{Sher69}.
 
E.R.~Apodaca 
embedded $I\times \mathcal C$ in $\mathbb R^3$
 so that $I\times \{s\}$, $s\in \mathcal C$,
are pairwise inequivalently 
embedded wild arcs,
each locally tame except at one point.
Also, identifying
$I\times \mathcal C$ 
with
$I\times \{ 0 \} \times \mathcal C
\subset I\times I\times \mathcal C $,
he extended his embedding to an embedding
$I\times I\times \mathcal C \to \mathbb R^3$ \cite{Apodaca}.

R.J.~Daverman constructed 
examples of $k$-cells $B$ in $\mathbb R^n$
($n\geqslant 4$, $3\leqslant k \leqslant n$)
such that 
each $2$-disk $D$ in $B$ (in $\partial B$ in case $k=n$)
 has non-simply connected complement in $\mathbb R^n$ and hence
  is wildly embedded in $\mathbb R^n$
 \cite{Daverman-absense}; this implies also a disjoint
embedding result.

Under certain restrictions,
if a space $Y$ contains uncountably many pairwise disjoint copies of a compactum $X$,
then $Y$ contains a copy of $X\times \mathcal C$
\cite{vanDouwen}, \cite{BEM}.
It would be interesting 
to specify this result
so that the copies $X\times \{s \}$
run over the prescribed set of
embedding types (up to ambient homeomorphism
of $Y$).

In a forthcoming paper,
instead
of Antoine's necklaces
we
will exploit Bing-Whitehead Cantor sets
which have simply-connected complements
but are nevertheless wild \cite[Theorem]{GO};
we will
construct
uncountable families of pairwise disjoint
wild $k$-disks in $\mathbb R^n$ such that
the complement
of each $k$-disk is simply connected;
the result generalizes \cite[Cor.~1]{GO}
and is announced 
in \cite[Thm.~1]{Frolkina-thesis}.

\section{Preliminaries}

A compactum is called
perfect if it has no isolated points.

By $\mathcal C$ we denote
the usual middle-thirds Cantor set  on $I=[0,1]$.
Any topological space
homeomorphic to $\mathcal C$
is called a Cantor set.
(These are exactly non-empty metric zero-dimensional perfect compacta, see \cite[Problem~1.3.F, p.~29]{Engelking}, 
\cite[Thm.~12.8]{Moise}.)

For a topological manifold-with-boundary
$M$,
denote by 
$\mathring M$ 
and $\partial M$
the interior and the boundary of $M$, correspondingly. In particular,
$\mathring I = (0,1)$.

For a subset $A$ of a topological space $X$,
its closure is denoted by
$\overline A$.

\begin{definition}\label{tame-C}
A zero-dimensional compact set $K\subset \mathbb R^n$ is called {\it tame\/}
if there exists
a homeomorphism $h$ of $\mathbb R^n$ 
onto itself such that 
$h(K)$ is a subset
of the $Ox_1$-axis of $\mathbb R^n$;
and it is called {\it wild\/} otherwise.

A zero-dimensional compact set $K\subset \mathbb R^n$ is {\it locally tame
at a point $x\in K$\/}
if
there is a neighborhood $U$ of $x$ in $\mathbb R^n$
such that
$K\cap\overline U$ is a tame (zero-dimensional) compactum, and {\it locally
wild at $x$\/} otherwise.

A zero-dimensional compact set $K\subset\mathbb  R^n$ is {\it everywhere wild\/}
if it is not locally tame at each of its points.
\end{definition}

A Cantor set $K$ in $\mathbb R^n$
is tame iff for some homeomorphism
$h:\mathbb R^n\cong \mathbb R^n$,
we have
$h(K) = \{ 0\}^{n-1} \times \mathcal C$
(refer e.g. 
to \cite[Prop.~6.1.17]{Bogachev}).

Each countable, compact subset of $\mathbb R^n$ 
is tame
\cite[1.1, 1.2]{Ivanov-diss},
\cite[Thm.~I.4.2]{Keldysh}.

We introduce the following

\begin{definition}\label{hyperplanar}
Call a subset
$X\subset\mathbb  R^n$
{\it locally hyperplanar
at a point $x\in X$\/}
if
there exist 
a neighborhood $U$ of $x$ in $\mathbb R^n$
and an embedding
$h:U\to \mathbb R^n$ 
such that 
$h (X\cap U) $
is contained in the standard
hyperplane
$\mathbb  R^{n-1}\times \{ 0\}$ of $\mathbb R^n$.
\end{definition}

A zero-dimensional compact subset of
a
hyperplane
$K\subset \mathbb R^{n-1} \times \{0\} \subset  \mathbb R^n$
is tame in $\mathbb R^n$, see
\cite[Thm.~2]{ZS}
(using Klee flattening theorem
\cite[Thm.~2.5.1]{Rushing}) or
\cite[Theorem~3]{Osborne1},
\cite[Cor.~1]{McMillan-Taming}.
Hence, in the next statement,
(a) implies (b);
the converse implication
is evident.

\begin{statement}\label{locally-wild-Cantor}
Let $n\geqslant 3$.
For any zero-dimensional compact set $K\subset\mathbb  R^n$
and a point $x\in K$ the following conditions are equivalent:

(a) $K$ is locally wild at $x$;

(b) $K$ is not locally hyperplanar at $x$.
\end{statement}

L.~Antoine in \cite{Antoine}
sketched and in
\cite[{\bf 78}, p.~91--92]{Antoine-diss}
explicitly constructed
a Cantor set 
in $\mathbb R^3$ which is now widely known as
Antoine's necklace. 
Antoine proved that
this Cantor set is wild
\cite[Part 2, Chap.~III]{Antoine-diss}.
Using ``feelers'' he 
showed that each zero-dimensional compactum in $\mathbb R^n$, $n\geqslant 2$, 
can be included in
a Jordan arc
 \cite[{\bf 72}, p.~82--84]{Antoine-diss} (this result was
announced by M.~Denjoy in \cite{D1}, \cite{D2}); 
thus Antoine described a 
first
example of a wild arc in $\mathbb R^3$
\cite[{\bf 83}, p.~97]{Antoine-diss}
(see also \cite[Thm.~18.8]{Moise}).
In $\mathbb R^n$, for each $n\geqslant 3$,
many other examples
of wild Cantor sets are now known.
We will use generalized Antoine's Necklaces 
as a base for our proofs (see Section
\ref{ABN}).

We would like to mention that besides
widely known Antoine's necklace,
there is an essentially different construction
of a wild Cantor set in $\mathbb R^3$ 
given 
by P.S.~Urysohn
\cite[n.40, p.~121--122]{Urysohn},
see also
\cite[p.~330--332]{Urysohn1}. 
According to P.S.~Alexandroff
\cite[p.~488]{Urysohn1},
P.S.~Urysohn wrote this paper in 1922--23.
Urysohn's set has infinite genus, in contrast to that of Antoine.
Urysohn writes with regret 
\cite[p.~332, footnote]{Urysohn1}
that 
Antoine's thesis came to his knowledge
only after he had almost finished his memoire.

\begin{definition}
\cite[Def.~1]{Keldysh62},
\cite[Def.~I.3.2]{Keldysh}
A set $K$ in an
$n$-dimensional topological manifold $M$
is called
{\it cellularly separated\/} (in $M$) if
for each open neighborhood $V$  of $K$
there exists 
a family $\{ u_\alpha \}$ of open subsets of $M$
such that
$K\subset 
\bigcup\limits_{\alpha }  u_\alpha 
\subset 
\bigcup\limits_{\alpha }  \overline u_\alpha  
\subset V$;
each $u_\alpha $ is homeomorphic to $\mathbb R^n$;
each
$\overline u_\alpha $
is a topological $n$-cell;
and 
$\overline u_\alpha  \cap \overline u_\beta = \emptyset $
for $\alpha \neq \beta $.
\end{definition}

The next proposition was
proved by L.V.~Keldysh
\cite[Thm.~1]{Keldysh62}, see also \cite[Thm.~I.4.2]{Keldysh}
(for $n=3$, the equivalency $(a)\Leftrightarrow (c)$ follows also from 
\cite[Thm.~1.1]{Bing-tame};
for arbitrary $n$, see also a later paper
\cite[Thm.~1]{Osborne1}).

\begin{statement}\label{cell-tame}
For each $n$ 
and each zero-dimensional
compact set
$K\subset \mathring I^n \subset \mathbb R^n$ the
following conditions are equivalent:

(a) $K$ is cellularly separated in $\mathring I^n$;

(b) there exists an isotopy $F_t$ of
$\mathbb R^n$ onto itself such that $F_0 = \id $;
 $F_t =\id $ on $\mathbb R^n - \mathring I^n$ for each $t\in I$; and 
$F_1(K) \subset Ox_1$;

(c) $K$ is tame in $\mathbb R^n$.
\end{statement}

The next Proposition
follows from the proof
of
\cite[Thm.~8]{Osborne1} (for $n=3$ 
see also \cite[Thm.~6.1]{Bing-tame}):

\begin{statement}\label{finite-sum}
If a Cantor set $K\subset\mathbb R^n$
is a union
of countably many 
tame compact sets $K_1, K_2,\ldots $,
then $K$ is tame.
\end{statement}

\begin{definition}\label{tame-P}
A subset $P\subset \mathbb R^n$
is called a {\it polyhedron\/}
if it is a union of a finite collection of simplices.
A compactum $X \subset\mathbb  R^n$  
homeomorphic to a polyhedron is
called
{\it tame\/}
if there exists a homeomorphism  
$h$ of $\mathbb R^n$ onto itself 
such that $h(X)$ is a polyhedron in $\mathbb R^n$;
and called {\it wild\/} otherwise.
\end{definition}

\begin{remark}
Any orientation-preserving self-homeomorphism of $\mathbb R^n$
is isotopic to the identity map.
This follows from the Alexander theorem
\cite[Thm.~V.3.1]{Keldysh}
together
with the Annulus Theorem
whose history is presented in
\cite{Cernavsky}.
This means that in Definitions
\ref{tame-C}
and \ref{tame-P} 
one can replace ``a homeomorphism of $\mathbb R^n$''
by ``an isotopy of $\mathbb R^n$''.
(For zero-dimensional compacta refer also to
Statement \ref{cell-tame}.)
\end{remark}

\begin{remark}
M.A.~Shtan'ko 
presented 
a general definition of tameness
for embedded compacta and studied
its properties
\cite{Shtanko69}, \cite{Shtanko}.
A brief but detailed overviews
analyzing its relations with other possible definitions can be found in
\cite{Me} and \cite{Cernavsky}.
\end{remark}

As we said, Antoine constructed
the first wild arc in $\mathbb R^3$.
In \cite{Antoine-173},
he announced
and in \cite{Antoine-FM}
described in detail 
a 2-sphere which contains 
the Antoine's necklace
(the description can also be found in
\cite{Alexander}
and 
in \cite[Thm.~18.7]{Moise}).
Antoine mentions that a homeomorphism 
between this 2-sphere 
and a standard 2-sphere 
in $\mathbb R^3$ can be extended
over their interiour domains, but cannot be extended over their
exterior domains, in constrast to the case of 
simple closed curves in plane;%
\footnote{
``It en r\'esulte qu'il existe dans $E_3$ une
surface hom\'eomorphe \`a une sph\`ere, la
correspondance s'\'etendant aux int\'erieurs 
de ces surfaces, mais pas aux ext\'erieurs... 
On sait que, au contraire, la correspondance entre deux courbes de Jordan planes peut s'\'etendre \`a la totalit\'e de leurs plans.'' \cite[p.~285]{Antoine-173}}
this homeomorphism cannot be extended
over neighborhoods of the 2-spheres.%
\footnote{``{\it ...qu'il existe, dans l'espace
\`a 3 dimensions, des surfaces hom\'eomorphes
(m\^eme simplement connexes), telles qu'aucune correspondance entre ces surfaces ne peut \^etre \'etendue \`a leurs voisinages.}''
\cite[p.~282--283]{Antoine-FM}}
A better-known 
example of 
a wild $2$-sphere in $\mathbb R^3$ 
is the 
horned sphere of J.W.~Alexander
\cite{Alexander-sphere}.

The feelers idea is very fruitful and widely 
used.
As said, it comes back to
\cite{Antoine-173},
\cite{Antoine-diss},
\cite{Antoine-FM},
\cite{Alexander},
\cite{Ivanov},
\cite{Blankinship},
\cite[proof of Theorem~3, Corollary~4]{Osborne2},
\cite[Thm.~1]{ZS}
and it can be stated as follows:

\begin{statement}\label{feelers-ext}
For $n\geqslant 3$,
any Cantor set $K$ which is a cellularly separated subset of $ \partial I^n$
and for any embedding
$f: K \to \mathbb R^n$,
there exists an embedding
$F:I^n \to \mathbb R^n$ such that
$F|_{K} = f$,
and $F$ is piecewise linear on
$I^n - K$.
\end{statement}
  
\begin{definition} 
Two subsets $X,Y\subset \mathbb R^n$ 
are {\it ambiently homeomorphic\/}
(or {\it equivalently embedded\/})
if there exists a homeomorphism $h$
of $\mathbb R^n$ onto itself such that
$h(X)=Y$.
\end{definition}

\begin{definition}
Let $X,Y\subset \mathbb R^n$.
We say that $X$ {\it can be ambiently embedded\/}
in $Y$ if there exists a homeomorphism $h$ of $\mathbb R^n$ onto itself such that $h(X)\subset Y$.
Two subsets $X,Y\subset\mathbb  R^n$ 
are called {\it ambiently comparable\/} if
at least one of them can be ambiently embedded into the other.
\end{definition}

\begin{remark}
Using Statements \ref{Sher-classif} and \ref{A-embedding}
and their higher-dimensional analogues
(Subsection \ref{ABN}),
one can easily construct for each $n\geqslant 3$
two wild Cantor sets $X,Y\subset \mathbb R^n$
such that $X$ ambiently embeds in $Y$,
and $Y$ ambiently embeds in $X$, but
$X$ and $Y$ are inequivalently embedded.
\end{remark}

\begin{remark}\label{obstr}
Suppose that $X,Y\subset \mathbb R^n$ and
$X$ ambiently embeds in $Y$. If $X$ is not locally hyperplanar at
a point
$x$, then $Y$
is not locally hyperplanar at the point
$h(x)$.
\end{remark}

\section{Reduction of Theorem~\ref{Thm2} to Theorems~\ref{MainLemma1}
and~\ref{reembedding}}\label{reduc}

By Theorem~\ref{reembedding}
we may assume that $X\subset I^{n-1}$
and $X\cap  \partial I^{n-1}
= \{ 0\}^{n-2} \times \mathcal C$,
where $\mathcal C$ is the middle-thirds Cantor set
on $I$.
Hence $
\{ 0\}^{n-2}
\times \mathcal C\times \mathcal C
\subset X\times \mathcal C\subset I^{n-1}\times I$.
Take an embedding
$$
\{ 0\}^{n-2}
\times \mathcal C\times \mathcal C
\cong \mathcal C\times\mathcal C
\stackrel{A}{\hookrightarrow }
\mathbb R^n,
$$
where 
the first identification is simply
$(0,\ldots , 0,x,y)\mapsto (x,y)$, and
$A$ is the embedding given by Theorem~\ref{MainLemma1}.
Extend it to
an embedding
$\tilde A : I^n \to \mathbb R^n$
which  is piecewise linear on
$ I^n - \{0\}^{n-2}\times\mathcal C\times\mathcal C$ (Statement~\ref{feelers-ext}).
Define the 
embedding $\mathbb A$ as the restriction
of $\tilde A$ on
$X\times \mathcal C$.
Let us show that $\mathbb A$ has the desired
properties.

Fix any $s\in\mathcal C$.
Evidently
$\mathbb A (X\times \{s\})
\supset
A(\mathcal C\times \{ s\})$, hence
 the set
$\mathbb A(X\times \{s\} )$ 
is not locally hyperplanar
at any point $x\in \mathbb A(
\{ 0 \}^{n-2}\times \mathcal C\times \{s\} )$
(see Statement~\ref{locally-wild-Cantor}
and Remark \ref{obstr}).
Since $\tilde A:I^{n}\to \mathbb R^n$
is piecewise linear on
$I^n- \{ 0\}^{n-2}\times \mathcal C\times \mathcal C$,
its restriction
on
$I^{n-1}\times \{s\} $ 
is piecewise linear on
$I^{n-1}\times \{s\} - \{ 0\}^{n-2}\times \mathcal C\times  \{ s\}$.
By \cite[Thm.~1.7.2]{Rushing}, the set
$\tilde A (I^{n-1}\times \mathcal \{ s\} )$
(and therefore
$\mathbb A (X\times \mathcal \{s\} )$)
is locally hyperplanar
at each point $x\in \mathbb A ( I^{n-1}\times \{s\} - \{ 0\}^{n-2}\times \mathcal C
\times \{ s\})$.
By Remark~\ref{obstr},
the sets 
$\mathbb A (X\times \{s\})$ for
$s\in\mathcal C$
are pairwise ambiently incomparable.

To prove {\it (i)},
note that
$\tilde A |_{X\times I}$
is the required isotopy.

It remains to prove {\it (ii)}.
Now $X$ is a polyhedron with $\dim X \geqslant 1$.
Suppose that
$\mathbb A (X\times \{s\})$ is tame 
for some $s\in \mathcal C$.
Hence
there exists
a polyhedron $P$ in $\mathbb R^n$
such that:
$P$ contains a Cantor set $K$ which
is everywhere wild as a subset of $\mathbb R^n$,
$\dim P \geqslant 1$,
and 
$P$ is embeddable in $\mathbb R^{n-1}$.
Represent $P$ as a finite union of simplices
$P=\bigcup \limits_{i=1}^{q} \sigma _{i}$.
(We do not assume that the simplices $\sigma _i$ are
pairwise disjoint.)
For each $i=1,\ldots , q$
we have $\dim \sigma _i \leqslant n-1$;
by \cite[Thm.~2]{ZS} the zero-dimensional 
compact set
$K\cap \sigma _i$ is tame in $\mathbb R^n$.
By
Statement \ref{finite-sum}
$K$ is tame, a contradiction.

Theorem~\ref{Thm2} is proved.

\section{Proof of Theorem \ref{MainLemma1}}\label{MainLemma1proof}

The idea is to construct ``the 
Cantor set of Antoine's necklaces
$\cup \mathcal A_s$'' and the embedding $A$
simultaneously.
We proceed in a countable number of steps.
In our construction of
$\cup \mathcal A_s$, we alternate steps of two types:
ramification and inserting simple chains.
The resulting set $\cup \mathcal A_s$
can be considered as a generalization of 
ramified Antoine's Necklaces, see
\cite[p.~383]{Ea}, 
\cite{Daverman-absense},
\cite{Wright1983}.

\subsection{On Antoine's Necklaces}\label{AN}

In 1920-21, L.~Antoine constructed and investigated
\cite{Antoine}, \cite[p.~91--105]{Antoine-diss}
 his famous wild Cantor set
now called Antoine's Necklace.

Let us recall his construction.

\begin{definition}\label{torus-r-R}
Let $\Pi $ be a plane in $R^3$.
Let $D\subset \Pi $ be a disk of radius $r$
with center $Q$, and $\ell \subset \Pi $ a straight line such that $d(Q,\ell ) = R > r$.
A {\it standard solid torus\/}
is the solid torus of revolution 
$T$ generated by
revolving $D$ in $R^3$ about $\ell $.
The {\it central circle\/} of $T$
is the circle generated by rotating the point $Q$. The {\it center\/} of $T$ is the center
of its central circle.
\end{definition}

Antoine takes a standard
solid torus $T$ in $\mathbb R^3$;
assume that the ratio $\frac{r}{R}$
is small enough
(where $r$ and $R$ are numbers described
in Definition \ref{torus-r-R}).
There exists an integer $k$ (sufficiently large 
in comparison with $\frac{r}{R}$) such that
a simple chain of cyclically linked $k$ tori, each 
{\it geometrically similar\/} to $T$, can be placed inside $T$,
so that 
their centers lie on the central circle of $T$
and form a regular convex $k$-gon.
Then, he applies a similarity transformation
to place a chain of $k$ tori in the interior of each torus of the previous level, and so on.
Since diameters of the tori tend 
to zero, the limit set
is a Cantor set.
Antoine showed
\cite[Part II, Chap. III]{Antoine-diss} that
this Cantor set is wild.
Other proofs can be found e.g. in
\cite[Prop.~9.5]{D-dec},
\cite[Section~18]{Moise},
\cite[\S IV.4]{Keldysh}.

To generalize
this construction, one
may allow the tori to be non-standard 
solid tori,
that is, not necessary tori of revolution;
also, one may vary the number of
tori on each stage, see \cite{Sher68}.
For us, it will be convenient
to 
restrict ourselves by the following 
definition.

\begin{definition}\label{A-def}
{\it A simple chain\/}
in a standard solid torus $T\subset \mathbb R^3$
is a finite family 
$T_1,\ldots , T_q$, $q\geqslant 3$, 
of pairwise disjoint 
congruent
standard solid tori such that

1) $T_1\cup\ldots\cup T_q\subset \mathring T$;

2) centers of $T_1,\ldots, T_q$
are subsequent vertices of a regular convex $q$-gon inscribed in the
central circle of $T$;

3) $T_i$ and $T_j$ are linked for
$|i - j |\equiv 1 \ \mod q$, and are not linked
otherwise;

4) for each $i$, the central circle of $T_i$
is zero-homotopic in $T$.

That is, the chain $T_1,\ldots , T_q$ 
looks like a usual ``necklace'' which
winds once around the hole of $T$;
no one of the tori $T_i$ embraces the hole of $T$.

{\it An Antoine's Necklace\/} is a Cantor set in $\mathbb R^3$
which can obtained as an
intersection $\mathcal A=\bigcap\limits_{i\geqslant  0}  M_i$,
where
each $M_i$ 
is the union of a finite number
of pairwise disjoint standard solid tori 
such that:

A1) $M_0\subset\mathbb  R^3$ 
is a standard solid torus;

A2) $M_{i+1} \subset M_i$ for each $i$;

A3) for each $i$ and each component $T$ of $M_i$,
the intersection
$M_{i+1} \cap T$ 
is a union of
solid tori,
which form a simple chain in $T$.

The sequence $\{M_i\}$ is called 
{\it a canonical defining sequence\/} for $\mathcal A$.

(Since $\dim \mathcal A =0$, conditions
A1)--A3) imply that diameters
of components of $M_i$ tend to zero as
$i\to \infty $. Conversely,
if a sequence $\{ M_i\} $ satisfies A1)--A3) 
and diameters
of components of $M_i$ tend to zero as
$i\to \infty $, then the intersection
$\bigcap\limits_{i\geqslant 0} M_i$ is
a Cantor set.)
\end{definition}

Each Antoine's necklace
is everywhere wild;
see
\cite[{\bf 79--82}, p.~93--96]{Antoine-diss}.

\begin{remark}
Many authors describe an Antoine's necklace
as obtained by first taking 
a standard solid torus, then placing
a simple chain of four solid tori in its interior,
then again placing a simple chain of four tori
in each of these four tori, and so on.
A simple chain is often defined ``by a picture''.
This is not the original Antoine's construction.
It is true that the diameters 
of the tori can be made arbitrary small
in case we use four {\it non-standard\/} tori at each stage (in fact, it suffices to take two tori \cite{Bing52}), 
but this is not evident since
the tori may become ``more and more curved''
on each stage.
For simplicity of exposition,
we prefer to use ``a large number'' of standard tori, as Antoine did.
(For an interesting related result, see \cite{Z}.)
\end{remark}

R.B.~Sher classified
Antoine's necklaces up to ambient homeomorphism
(in fact, he considered
more general Antoine-type sets
than those described in Definition~\ref{A-def}):

\begin{statement}\label{Sher-classif}
\cite[Thm.~2]{Sher68}
Let $\mathcal A$, $\mathcal B$ be two Antoine's necklaces in $\mathbb R^3$
with canonical defining sequences $\{ M_i\}$, $\{ N_i\} $. 
Then $\mathcal A$ and $\mathcal B$ are ambiently homeomorphic 
if and only if
there exists a homeomorphism
$h$ of $\mathbb R^3$ onto itself such that
$h(M_i)=N_i$ for each $i\in\mathbb N$.
\end{statement}

\begin{remark}\label{Sher-weak}
Sher's Theorem implies that ambiently 
homeomorphic necklaces 
have the same number of components on each level
$i$.
This statement is weaker than Sher's Theorem itself
(see an example in
\cite[Theorem~2]{GRZ2005}), and its proof
can be obtained directly;
see \cite[Theorem 4.6]{Wright-rigid}.
\end{remark}

\begin{statement}\label{A-embedding}
Let $\mathcal A,\mathcal B \subset \mathbb R^3$ be two Antoine's necklaces
with canonical defining sequences 
$\{ M_i\} $ and $\{ N_i\}$ correspondingly.
Then,
 $\mathcal A$ can be ambiently embedded in $\mathcal B$
 if and only if there exist
 an integer $k\geqslant 0$,
 a component $T$ of $N_k$,
 and a self-homeomorphism
$h$ of $\mathbb R^3$ such that
$h(M_i) = T\cap N_{i+k}$ for each~$i\in \mathbb N$.
\end{statement}

Proof.
Let $h$ be a self-homeomorphism of $\mathbb R^3$ onto itself such that $h(\mathcal A)\subset \mathcal B$.
Theorem 4.6 of \cite{Wright-rigid}
implies that
for some component $T$ of $N_k$,
$h$ maps $\mathcal A$ [homeomorphically] onto
$T\cap \mathcal B$; now it remains to apply Sher's Theorem. The converse statement is evident.

\subsection{Notation and Tools needed for
the proof of Theorem~\ref{MainLemma1} }\label{Lemma-Notation}

1)
The usual ``middle-thirds''
Cantor set $\mathcal C$ in $[0,1]\subset \mathbb R^1$
is defined as 
follows.

Take $K_0 = I=[0,1]$.
Put $K_1 = \Delta _0\cup \Delta _1$, where
$\Delta _ 0 = [0,\frac13]$
and $\Delta _1 = [\frac23 , 1]$.
Suppose that we have constructed
a
family $ S (I;N)$ of
$2^N$ segments 
$\Delta _{i_1 i_2 \ldots i_N}$, where
$(i_1,\ldots ,i_N) \in \{ 0,1\}^N$; put
$K_N = \cup \{ J \ | \ J \in  S(I,N) \}$.
Divide each
$\Delta _{i_1 i_2 \ldots i_N}
\in  S (I;N)$
into 3 equal parts; define
$\Delta _{i_1 i_2 \ldots i_N 0}$
and $\Delta _{i_1 i_2 \ldots i_N 1}$
to be its first and third part,
 correspondingly
(going from left to right).
We thus obtain the family $ S (I;N+1)$
of
$2^{N+1}$ segments 
$\Delta _{i_1 i_2 \ldots i_N i_{N+1}}$,
where $(i_1,\ldots ,i_{N+1}) \in \{0,1\}^{N+1}$.
Define
$
K_{N+1}
= \cup \{ J \ | \ J \in  S (I;N+1) \} 
$, and so on.
Finally, put
$\mathcal C = \cap_{N=0}^\infty K_N$.

2)
For any
segment $L=[a,b]$ 
and any
  given integer $N\geqslant 0$,
define the family of $2^N$ segments on $L$
by the formula
\begin{equation}
\label{S}
 S (L;N) 
= \{ H_L(J) \ |  \ 
J \in  S (I;N)
\}
\end{equation}
where
$H_L:I\to L$
is the homeomorphism
defined by
$t\mapsto a(1-t) + bt $, $t\in I$.

3) Let $\{ L_i = [a_i,b_i]\}$
be a finite family of pairwise disjoint 
subsegments of $[0,1]$.
We say that this system is {\it regularly indexed}
or {\it regularly ordered}
if $i<j$ implies $b_i < a_j$.

4) For each
 non-negative integer $s$,
 each array
$(a_{i_1}, a_{i_1i_2},\ldots , 
a_{i_1\ldots i_s})$ of positive integers, and
each array
$$
(k_{i_1} ,
k_{i_1 i_2},\ldots ,
k_{i_1 i_2 \ldots i_s })
\in 
\{ 1,\ldots , 4^{a_{i_1}} \}\times
 \{ 1,\ldots , 4^{a_{i_1 i_2}} \}\times
\ldots 
\times
 \{ 1,\ldots , 4^{a_{i_1 i_2\ldots i_s}} \}
$$
define a segment 
$J\langle 
k_{i_1}, k_{i_1 i_2},\ldots , k_{i_1 i_2\ldots i_s}
\rangle $ inductively as follows.

For $s=0$, the unique segment is $J=I=[0,1]$.
Suppose that 
we have defined all segments
$J\langle 
k_{i_1}, k_{i_1 i_2},\ldots , k_{i_1 i_2\ldots i_{q-1}}
\rangle $ for each $q\leqslant s$.
Then, to obtain 
$J\langle 
k_{i_1}, k_{i_1 i_2},\ldots , k_{i_1 i_2\ldots i_{s}}
\rangle $,
consider   the segment
$J\langle 
k_{i_1}, k_{i_1 i_2},\ldots , k_{i_1 i_2\ldots i_{s-1}}
\rangle .$
Enumerate 
$4^{a_{i_1 i_2\ldots i_s}}$
segments of the family
$ S (J\langle 
k_{i_1}, k_{i_1 i_2},\ldots , k_{i_1 i_2\ldots i_{s-1}}
\rangle ,
2a_{i_1 i_2\ldots i_s}) $
(see formula (\ref{S}))
in regular order.
Define 
$J\langle 
k_{i_1}, k_{i_1 i_2},\ldots , k_{i_1 i_2\ldots i_s}
\rangle $
to be the 
$k_{i_1 i_2\ldots i_s}$-th 
segment of this family.

5) 
Let $G$ be a family of subsets of $\mathbb R^n$.
The number of its elements is denoted by $\# G$;
 the mesh of the family is 
$\| G \| = \sup  \{ \diam A \ | \ A\in G \}.$
By $\widetilde{G}$ we denote the union
$\cup \{ g \ | \ g\in G \}=
\{ x: \ x\in g \ \text{ for some }\ g\in G \}$.
For two families  $G, G'$ of subsets of $\mathbb R^n$, the inequality
$G\leqslant G'$ means that $G$ is a refinement of $G'$.

Next proposition follows from
 \cite[Theorem~12.7]{Moise}:

\begin{statement}\label{Moise}
Suppose that 
for each integer $k\geqslant 0$ we have 
families $G_k$ and $G'_k$
of subsets
of $\mathbb R^n$ and $\mathbb R^m$ correspondingly
such that

(a) both $G_k$ and $G'_k$ are finite families
of pairwise disjoint compacta;

(b) $G_{k+1} \leqslant G_k$
and $G'_{k+1} \leqslant G'_k$;

(c) $\# G_k = \# G'_k$;

(d) $\| G_k \| \to 0$ and
$\| G'_k \| \to 0$
as $k\to \infty $.

Moreover, suppose that for each $k\geqslant 0$ we are given a bijection $f_k : G_k \to G'_k$ such that
the triple of conditions
$g_k\in G_k$, $g_{k+1} \in G_{k+1}$,
$g_{k+1}\subset g_k$
implies $f_{k+1}(g_{k+1}) \subset f_k (g_k)$.

Then 
each point $P\in \bigcap \limits_{k\geqslant 0} \widetilde{G_k} $
has a unique representation of the form
$\{ P \} =\bigcap \limits_{k\geqslant 0} g_{k,P}$, where
$g_{k,P} \in G_k$ for each $k\geqslant 0$;
the formula
$\{ f(P)\} = \bigcap \limits _{k\geqslant 0} f_k (g_{k,P})$
provides a well-defined map $f$
of $\bigcap\limits _{k\geqslant 0} \widetilde{G_k}$ to 
$\bigcap \limits_{k\geqslant 0} \widetilde{G'_k}$;
and $f$ is a homeomorphism onto.
\end{statement}

Below we will need
the 2-fold Ramification procedure.

\begin{definition}
Suppose we are given a standard solid torus 
$T$ in $R^3$; $T$ is generated by revolving
a disk $D$ about an axis $\ell $
(Definition \ref{torus-r-R}).
Choose two disjoint subdisks $D_1$, $D_2 $ of
$\mathring D$
and take the
standard solid tori $T_1$, $T_2$
obtained by revolving $D_1$, $D_2$
about $\ell $.
We call
$T_1\cup T_2$
a {\it 2-fold ramification} of $T$.
\end{definition}

The following fact is elementary
(the number of links is even for simplicity;
we have in mind that the central circle
of every second link lie in the same plane
which contains the central circle of $T$,
and it is perpendicular to 
planes containing central circles of other links):
 
\begin{statement}\label{small-tori}
Let $T$ be a solid torus of revolution.
For each $\varepsilon >0$
there exists an integer $N$
such that for each $k\geqslant N$
one can construct in $T$ a simple
chain of $2k$ congruent standard solid tori
whose diameters
do not exceed $\varepsilon $.
\end{statement}

\subsection{Proof of Theorem \ref{MainLemma1} for $n=3$.}

The desired embedding $A$
will be constructed using
Statement \ref{Moise}.
We will
inductively define families $G_k$ and $G'_k$ 
of
subsets of $\mathbb R^2$ and $\mathbb R^3$ correspondingly, together 
with bijections $f_k$, so that
$\bigcap \limits_k \widetilde{G_k}  = \mathcal C \times \mathcal C$,
and $\bigcap \limits_k \widetilde{G_k'}$
is a Cantor set of pairwise
ambiently incomparable Antoine's necklaces.
For even and for odd values of $k$
the constructions are different;
we use therefore indices of two types:
the first ones are usual; and the second ones
are in square brackets.
The whole array of indices is written in angle brackets.

In the description, we will use
a special sequence
$a_0, a_1, a_{00}, a_{01}, a_{10}, a_{11}, \ldots $.
We assume that its members are
pairwise distinct positive integers.
We also assume that it increases 
sufficiently fast;
exact meaning of this will be clarified below.

Let us describe several initial steps for clarity.

\underline{Step 0.}
The family $G_0 $ consists of one element
$ K =  I\times I$.
The family $G'_0 $ contains one element $T$
which is a standard solid torus
in  $\mathbb  R^3 $.
The bijection $f_0: G_0\to G'_0$ is the obvious one.

\underline{Step 1.}
The family $G_1$ contains 2 elements:
$K\langle [0]\rangle  = I\times \Delta _0$
and
$K\langle [1]\rangle = I\times \Delta _1$.
The family $G'_1$ contains 2 elements $T\langle [0]\rangle $ and $T\langle [1]\rangle $ which are obtained from $T$ by the 2-ramification procedure.
Define the bijection
$f_1: G_1\to G'_1$ by 
$K\langle [i]\rangle \mapsto T\langle [i]\rangle $, $i=0,1$.

\underline{Step 2.}
Let $G_2$ 
contain elements of 2 types.
\\
Type 1: 
$K\langle [0];k_0\rangle  = 
J\langle k_0\rangle\times \Delta _0  $,
where
for 
$k_0=1,\ldots , 4^{a_0}$ the segments
$J\langle k_0\rangle  $ run over the family
$ S(I; 2a_0)$, see the 
formula $(\ref{S})$ above.
\\
Type 2:
$K\langle [1];k_1\rangle  = J\langle k_1\rangle \times \Delta _1 $,
where
for 
$k_1=1,\ldots , 4^{a_1}$ the segments
$J\langle k_1\rangle  $ run over the family
$ S(I; 2a_1)$.

Also, elements of $G'_2$ are of two types.
\\
Type 1:
solid tori
$T\langle [0];k_0\rangle $,
$k_0=1,\ldots , 4^{a_0}$, which form a
simple chain
 in $T\langle [0]\rangle $.
\\
Type 2:
solid tori
$T\langle [1];k_1\rangle $,
$k_1=1,\ldots , 4^{a_1}$,
which form a
simple chain 
in $T\langle [1]\rangle $.

We can assume that $a_0$ and $a_1$ are
sufficiently large,
so that
$\| G'_2 \| < \frac12$
(see Statement \ref{small-tori}).

The
bijection $f_2:G_2\to G'_2$
is defined by
$ K\langle [i];k_i\rangle  \mapsto T\langle [i];k_i\rangle $
for each $i=0,1$; $k_i=1,\ldots , 4^{a_i}$.

Let us now describe the general rule.

\underline{Step 2s+1.}
On Step $2s$ we constructed
the family $G_{2s}$. Its elements are of the form
$$
K\langle
[i_1], k_{i_1}, [i_2], k_{i_1i_2},
 \ldots ,
 [i_s],
k_{i_1 \ldots i_{s}}
\rangle 
=
J \langle
k_{i_1}, k_{i_1i_2},\ldots ,
k_{i_1\ldots i_{s}}
\rangle 
\times\Delta _{i_1\ldots i_s}
,
$$
where  each $i_j  \in \{ 0,1\}$, each
$k_{i_1\ldots i_j} \in \{ 1,\ldots , 4^{a_{i_1\ldots i_j}}\}$,
and 
the segments
$J \langle
k_{i_1}, k_{i_1i_2},\ldots ,
k_{i_1\ldots i_{s}}
\rangle $
are defined as in Section
\ref{Lemma-Notation}.

Define the family 
$G_{2s+1}$. Its elements are
$$
K\langle
[i_1], k_{i_1}, [i_2], k_{i_1i_2},
 \ldots ,
 [i_s],
k_{i_1 \ldots i_{s}}, [i_{s+1}]
\rangle 
=
J \langle
k_{i_1}, k_{i_1i_2},\ldots ,
k_{i_1\ldots i_{s}}
\rangle 
\times \Delta _{i_1 \ldots i_s i_{s+1}} 
,
$$
where each $i_j \in \{ 0,1\}$ and each
$k_{i_1\ldots i_j} \in \{ 1,\ldots , 4^{a_{i_1\ldots i_j}}\}$.

To obtain the family $G'_{2s+1}$,
apply the 2-ramification procedure
to each
${T\langle
[i_1], k_{i_1}, [i_2], k_{i_1i_2},
 \ldots ,
 [i_s],
k_{i_1 \ldots i_{s}}
\rangle }$;
denote the two resulting tori by
$$
T\langle
[i_1], k_{i_1}, [i_2], k_{i_1i_2},
 \ldots ,
 [i_s],
k_{i_1 \ldots i_{s}}, [0]
\rangle 
$$
and
$$
T\langle
[i_1], k_{i_1}, [i_2], k_{i_1i_2},
 \ldots ,
 [i_s],
k_{i_1 \ldots i_{s}}, [1]
\rangle .
$$
The family $G'_{2s+1}$
consists of all
solid tori 
$$
T\langle
[i_1], k_{i_1}, [i_2], k_{i_1i_2},
 \ldots ,
 [i_s],
k_{i_1 \ldots i_{s}}, [i_{s+1}]
\rangle ,
$$
where indices run over the same domains as for $G_{2s+1}$.
Note that the inequality $\| G'_{2s} \| < \frac{1}{2s}$
implies $\| G'_{2s+1} \| < \frac{1}{2s}$.

The bijection
$f_{2s+1}: G_{2s+1} \to G'_{2s+1}$
maps
$$
K\langle
[i_1], k_{i_1}, [i_2], 
 \ldots ,
 [i_s],
k_{i_1 \ldots i_{s}}, [i_{s+1}]
\rangle  
\mapsto
T\langle
[i_1], k_{i_1}, [i_2], 
 \ldots ,
 [i_s],
k_{i_1 \ldots i_{s}}, [i_{s+1}]
\rangle  .
$$

\underline{Step 2s+2.}
The family $G_{2s+2}$
consists of the elements
$$
K\langle
[i_1], k_{i_1}, [i_2], k_{i_1i_2},
 \ldots ,
[i_{s+1}],
k_{i_1 \ldots i_{s+1}}
\rangle 
=
J \langle
k_{i_1}, k_{i_1 i_2 },
\ldots ,
k_{i_1\ldots i_{s+1}}
\rangle
\times
\Delta _{i_1 \ldots i_{s+1} } 
 ,
$$
where each $i_j \in \{ 0 ,1\} $, and each
$k_{i_1\ldots i_j} \in \{ 1,\ldots ,
4^{a_{i_1\ldots i_j}} \}$.

Let us describe the family $G'_{2s+2}$.
For
any 
fixed values of
$[i_1]$, $k_{i_1}$, $[i_2]$, $k_{i_1i_2},
 \ldots ,
[i_{s+1}]$,
where each $i_j \in \{0,1\}$ and each
$k_{i_1\ldots i_j} \in \{1,\ldots , 
4^{a_{i_1 \ldots i_j}} \}$,
take
a simple chain
of
$4^{a_{i_1 \ldots i_{s+1}}}$
standard solid tori
  in
$T\langle
[i_1], k_{i_1}, [i_2], k_{i_1i_2},
 \ldots ,
 [i_{s+1}]
\rangle 
$.
Denote these tori by
$T\langle
[i_1], k_{i_1}, [i_2], k_{i_1i_2},
 \ldots ,
[i_{s+1}],
k_{i_1 \ldots i_{s+1}}
\rangle 
$,
where
$ k_{i_1 \ldots i_{s+1}} \in \{ 1,\ldots ,4^{a_{i_1 \ldots i_{s+1}}} \}$.
All tori obtained in this way form
the family $G'_{2s+2}$.

We can assume that the numbers $a_{i_1\ldots i_{s+1}}$ are sufficiently large,
so that
$\| G'_{2s+2} \| < \frac{1}{2s+2}$
(see Statement \ref{small-tori}).

The bijection
$f_{2s+2}: G_{2s+2} \to G'_{2s+2}$
is defined by
$$
K\langle
[i_1], k_{i_1}, [i_2], 
\ldots ,
 [i_{s+1}], k_{i_1 \ldots i_{s+1}}
\rangle  
\mapsto
T\langle
[i_1], k_{i_1}, [i_2], 
 \ldots ,
 [i_{s+1}], k_{i_1 \ldots i_{s+1}}
\rangle  .
$$

The process continues ad infinitum.

By construction,
$\| G'_k \| \to 0 $
as $k\to \infty $.
Hence 
$K' = \bigcap\limits_k \widetilde{G'_k}$ 
is a Cantor set.
The intersection
$\bigcap \limits_k \widetilde{G_k}$ 
is exactly $\mathcal C \times \mathcal C$.

Statement \ref{Moise}
provides 
the homeomorphism $\mathcal C \times \mathcal C \cong K' \subset \mathbb R^3$; thus we obtain the
desired embedding 
$A:\mathcal C \times \mathcal C\to \mathbb R^3$.

In fact, each
$s\in \mathcal  C$
is the intersection
of a (unique) system
$$\Delta _{i_1} \supset \Delta _{i_1 i_2} \supset \ldots $$
where each $i_j \in \{ 0,1\}$.
The image $A(\mathcal C \times \{ s \} ) $
is easily seen to be the intersection
of the corresponding decreasing sequence of sets
$$
T\langle [i_1]\rangle 
\supset
\bigcup\limits_{k_{i_1} = 1,\ldots , 4^{a_{i_1}}}
T\langle [i_1], k_{i_1}\rangle 
\supset
$$
$$
\supset
\bigcup\limits_{k_{i_1} = 1,\ldots , 4^{a_{i_1}}}
T\langle [i_1], k_{i_1}, [i_2] \rangle 
\supset
\bigcup_{\substack{k_{i_1} = 1,\ldots , 4^{a_{i_1}}
\\
k_{i_1 i_2 } = 1,\ldots , 4^{a_{i_1 i_2}}}
}
T\langle [i_1], k_{i_1}, [i_2], k_{i_1 i_2} \rangle 
\supset
\ldots  
$$
that
is,
$A(\mathcal C \times \{ s \} ) = : \mathcal A_s$ is an Antoine's necklace in $\mathbb R^3$.

Recall that the numbers 
$a_{i_1\ldots i_k}$ are pairwise different.
By Statement \ref{A-embedding} the
necklaces $\mathcal A_s$ are pairwise
ambiently incomparable.

\subsection{Proof of Theorem \ref{MainLemma1} for $n\geqslant 4$ (sketch).}\label{ABN}

This case has minor difference with
the case of $n=3$.
Higher-dimensional generalizations
of Antoine's construction
were suggested
by 
A.A.~Ivanov \cite{Ivanov},
\cite{Ivanov-diss} and by
W.A.~Blankinship \cite{Blankinship};
they both used
$n$-dimensional solid tori
$D^2\times (S^1)^{n-2}$
in $\mathbb R^n$ and
showed that 
the complements of these ``generalized
Antoine's necklaces'' have
non-trivial fundamental groups.
Analogous to Antoine's construction,
Blankinship takes
the same number $k$ of smaller tori inside
each torus on each level, on each stage 
they are placed likewise on the first stage.
A broader class of sets can be obtained
by varying numbers of tori constructed 
on each level in each component.
(Ivanov's scheme allows this.)
W.T.~Eaton gives a brief and very clear description of such
generalized Antoine-type sets 
\cite[p.~380--381]{Ea}
and proves that they are wild \cite[p.~381--383]{Ea}.
(Suchlike sets were described 
and studied also in \cite{Osborne3}.)
Each of the sets described by Eaton
is everywhere wild; further, for these sets
we have
a weaker but satisfactory analogue of Statement \ref{A-embedding},
see \cite[Theorems~4.6, 5.1]{Wright-rigid}.

Having said all this, we see that
a detailed proof for $n\geqslant 4$
is a hand-by-hand repetition of 
the case $n=3$. It is therefore omitted.

\section{Proof of Theorem~\ref{reembedding}}

Let us introduce some additional notation.
An {\it $n$-cube neighborhood} $N$
of a  point
$p \in \mathbb R^n$ 
is a set of the form
$N = 
[a_1 , b_1] \times 
\ldots \times 
[a_n , b_n ]
$
such that 
each $b_i-a_i>0$ and $p\in \mathring N$.

For an $n$-cube neighborhood $N$ 
of both points $p, p'$
define a homeomorphism $H_{N, p,p'}:
N\cong N$
as follows.
Every point in $N$, other than $p$,
lies on a unique line segment
$pq$ where $q\in \partial N$.
The homeomorphism
$H_{N,p,p'}$ takes $p$ to $p'$;
and it maps the line segment $pq$
onto $p'q$ linearly.
Note that $H_{N,p,p'}$ is the identity on $\partial N$. 
Extending $H_{N,p,p'}$ by the identity map,
we get a homeomorphism $\widehat {H}_{N,p,p'} : \mathbb R^n \cong \mathbb R^n$.

For $n\geqslant 2$,
let
$\pi  : \mathbb R^n  \to \mathbb R^{n-1}$
and 
$\lambda  : \mathbb R^n  \to \mathbb R$
be the projection maps
given by
$\pi (x_1, \ldots , x_{n-1} , x_n) =
(x_1, \ldots , x_{n-1} )
$
and
$\lambda (x_1, \ldots , x_{n-1} , x_n)  = x_n $.

Proof of the next Lemma is straightforward.

\begin{lemma}\label{1}
In $\mathbb R^n$, $n\geqslant 2$,
let $N$ be an $n$-cube neighborhood
of two points, $p$ and $q$ with 
$\lambda (p) > \lambda (q)$.
We suppose that $x\in N$ and $\lambda (x) \geqslant \lambda (p)$.
Then 
$\lambda (H_{N,p,q} (x)) \geqslant 
\lambda (q)$; and equality holds
if and only if $x=p$.
\end{lemma}

\begin{lemma}\label{2}
Let $X\subset \mathbb R^n$ be a non-empty
compact set, $n\geqslant 2$.
Then there exists a homeomorphism $h:\mathbb R^n\cong \mathbb R^n$
such that
$h(X) \subset I^n$,
and
$h(X) \cap \partial I^n$
is a single point in 
$\mathring I^{n-1} \times \{0\}$.
\end{lemma}

Proof.
Since $X$ is bounded,
a linear homeomorphism
$L:\mathbb R^n\cong \mathbb R^n$
takes $X$ into $\mathring I^n$.
Let $\delta = \inf \limits _{p\in L(X)} \lambda (L(X))$. Recall that $L(X)$ is compact; 
hence $\delta >0$,
and for some point $p\in L(X)$
we have $\lambda (p) = \delta $.
Let $N = I^{n-1}\times [-1,1]$.
Take $p' \in \mathring I^{n-1} \times \{0\}$.
By Lemma \ref{1},
$\widehat {H}_{N,p,p'}\circ L$
is the desired homeomorphism.

\begin{statement}\label{feelers}
Let $X\subset I^n$ be a non-empty perfect set,
$n\geqslant 2$.
Suppose that
$X\cap \partial I^n = 
\{ p_1 , \ldots , p_k \}
\subset \mathring I^{n-1} \times \{ 0\}$.
For each $i$, let $N_i$ an $n$-cube neighborhood
of $p_i$ in $\mathbb R^n$
such that $\pi (N_i) \subset \pi (I^n)$,
and $N_i\cap N_j =\emptyset $ if $i\neq j$.
Then there exists a 
homeomorphism
$h:\mathbb R^n\cong \mathbb R^n$ such that

1) outside $\bigcup\limits _{i=1}^k N_i$,
$h$ is the identity map;

2) $h(X) \subset I^n$;

3) $h(X) \cap \partial I^n =
\{ p_1, q_1, p_2,q_2 , \ldots , p_k,q_k\}$,
where $q_i\in\pi(\mathring N_i)\times \{ 0\}$ and
$p_i\neq q_i$ for each $i$.
\end{statement}

Proof.
Let is fix an $i$ and construct the homeomorphism in $N_i$. For simplicity,
we omit the index $i$ and write
$N=N_i$ and $p=p_i$.
Since $X$ is perfect,
there exists a point $r\in X\cap \mathring N\cap \mathring I^n$. Consider two cases.
\\
{\it Case 1.} A point $r$ can be chosen so that
$\pi (r) \neq \pi (p)$.
Let $q = (\pi (r) , 0)$. Let $M\subset  N$
be an $n$-cube neighborhood
of both $r$ and $q$ which
does not contain $p$.
The set $X\cap M$ is compact and non-empty,
hence $\delta  = \inf \limits_{x\in X\cap M}
\lambda (x)$ is positive, and for some
$s\in X\cap M$
we have $\lambda (s) = \delta $.

{\it Subcase 1.1.}
Suppose 
that this point $s$ 
can be chosen in $X\cap \mathring  M$.
Then, $\widehat {H}_{M,s,q}$ is the desired homeomorphism.

{\it Subcase 1.2.}
Suppose 
that such a point $s$ 
cannot be chosen in $X\cap \mathring  M$.
Take any point $s'\in \mathring M$ with
$\lambda (s') = \delta $.
We come to Subcase 1.1
replacing $X$ by $\widehat {H}_{M,r,s'} (X)$
(the number $\delta $ may have changed).
\\
{\it Case 2.} 
For each $r\in X\cap \mathring N\cap \mathring I^n$, we have
$\pi (r) = \pi (p)$.
Let $M$
be an $n$-cube neighborhood of $r$
such that $M\subset \mathring N\cap \mathring I^n$.
Take any point $r' \in \mathring M$ with 
$\pi (r' ) \neq \pi (p)$.
Replace $X$ by
$\widehat {H}_{M,r,r'} (X)$; Case 1 now applies.

\begin{statement}\label{limit-emb}
Let 
$X\subset \mathbb R^n$ be a non-empty
perfect compact set, $n\geqslant 2$.
There exists a homeomorphism
$h:\mathbb R^n \cong \mathbb R^n$ such that

1) $h(X)\subset I^n$;

2) $h(X) \cap \partial I^n\subset 
\mathring I ^{n-1} \times \{ 0 \}$;

3) $h(X) \cap \partial I^n$ 
is a Cantor set 
cellularly separated as a subset of
$\mathring I ^{n-1} \times \{ 0 \}$.
\end{statement}

Proof.
The homeomorphism $h$ is constructed as a limit
of a sequence of homeomorphisms of
$ \mathbb R^n$.
This sequence is defined inductively.

By Lemma~\ref{2} we may assume 
that
$X \subset I^n$
and
$X \cap \partial I^n$ is a single point 
in $\mathring I ^{n-1} \times \{ 0 \}$;
call it $p$.
Let $N$ be an $n$-cube neighborhood of $p$
such that $\pi (N) \subset \mathring I ^{n-1} $. 
Apply Statement~\ref{feelers} 
to $X$ and $N$;
we obtain a homeomorphism $h_0: \mathbb R^n
\cong \mathbb R^n $ such that 
$h_0=\id $ outside $N$;
$h_0(X) \subset I^n$; and
$h_0(X) \cap \partial I^n  = \{ p_0, p_1\}
\subset \pi (\mathring N) \times \{ 0\} 
$, where $p_1\neq p_0 = p$.

Let $N_0$, $N_1$ be disjoint 
$n$-cube neighborhoods
of points $p_0$, $p_1$ correspondingly
with
$N_i \subset \mathring N$.
Apply Statement~\ref{feelers} 
to $h_0(X)$, $N_0$, $N_1$;
we obtain a homeomorphism
 $h_1:\mathbb R^n\cong \mathbb R^n$
 with 
 $h_1=\id $ outside $N_0\sqcup N_1$;
$h_1(h_0(X)) \subset I^n$; and
$h_1(h_0(X)) \cap \partial I^n 
= \{ p_{00}, p_{01}, p_{10}, p_{11} \}$, where 
$p_{00} = p_0 \neq p_{01} \in \pi (\mathring N_0 )\times \{ 0\}$, 
$p_{10} = p_1 \neq p_{11} \in \pi (\mathring N_1)
\times \{ 0\}$.

For each $(i_1,i_2)\in \{0,1\}^2$,
let
$N_{i _1i_2}$ be an $n$-cube neighborhood
of $p_{i_1i_2}$ such that 
$N_{i _1 0}\cup N_{i_1 1} \subset \mathring N_{i_1}$, and
$N_{i _1 0}\cap N_{i_1 1} = \emptyset $.
Apply Statement~\ref{feelers} 
to $h_1(h_0(X))$ and four neighborhoods $N_{i_1i_2}$; get a homeomorphism $h_2$.

Continuing in this way, we obtain for each $k\geqslant 1$
and each $(i_1,\ldots , i_k) \in \{ 0,1\}^k $
a homeomorphism
$h_k : \mathbb R^n\cong\mathbb R^n$,  
and a family
of $2^k$ pairwise disjoint
$n$-cube neighborhoods $N_{i_1\ldots i_k}$
 of
points
$p_{i_1\ldots i_k}$, such that
\begin{enumerate}
\item
$h_k = \id $ outside $\bigcup\limits_{(i_1\ldots i_k)\in \{0,1\}^k} N_{i_1\ldots i_k}$;
\item
$N_{i_1 \ldots i_k 0} \cup
N_{i_1 \ldots i_k 1 }\subset \mathring N_{i_1\ldots i_{k}};
$
\item
$h_{k-1} \circ h_{k-2}\circ \ldots \circ h_0 (X)
\subset I^n$;
\item
$h_{k-1} \circ h_{k-2}\circ \ldots \circ h_0 (X) \cap \partial I^n = \{ p_{i_1\ldots i_k} \ |\ (i_1,\ldots , i_k) \in \{ 0,1\}^k \} $;
\item
$ p_{i_1\ldots i_k} \in \pi(\mathring N_{i_1\ldots i_k})\times \{ 0 \}$.
\end{enumerate}
Note that for each $k$
we can choose the $n$-cubes $N_{i_1\ldots i_k}$
sufficiently small; how small depends
on the homeomorphism $h_{k-1}\circ h_{k-2}\circ\ldots\circ h_0$;
thus 
the sequence 
$\{ h_k \circ h_{k-1}\circ \ldots \circ h_0  \}$
can be constructed so that it
converges to a homeomorphism
$h:\mathbb R^n\cong \mathbb R^n$
\cite[Thm. 5.1]{Chapman}.
We have
$$
h(X) \cap \partial I^n = 
\bigcap\limits_{k\geqslant 1} \bigcup\limits_{(i_1,\ldots , i_k)\in \{0,1\}^k} N_{i_1\ldots i_k} \subset
\mathring I^{n-1}\times \{0\},
$$
this is a Cantor set cellularly separated
in $\mathring I^{n-1}\times \{0\}$.

\noindent
{\sc Proof of Theorem~\ref{reembedding}.}
We may assume that $X$
satisfies the conclusion 
of Statement \ref{limit-emb}.
That is, 
$X \subset I^n$ and
$K = X \cap \partial I^n$ 
is a cellularly separated
Cantor set
in 
$\mathring I^{n-1} \times \{ 0\}$.
There exists a homeomorphism
$G: \partial I^n \cong\partial I^n$
such that
$G(K) = \{ 0\}^{n-1} \times \mathcal C$
(see Statement~\ref{cell-tame}).
Extend $G$ to a homeomorphism 
$F : \mathbb R^n\cong \mathbb R^n$;
we get the desired map.

\medskip

{\bf Acknowledgement.}
The author is very much indebted to the referee for important comments on the first version of this paper.
Section 5 in its present form is
essentially due to the referee.
In my preprint,
statements and proofs were complicated.
I used ``topological'' feelers;
the referee suggested
``piecewise linear'' ones;
this allows to ``control''
their lowest points (Lemma~\ref{1}),
and the whole process becomes more transparent.
The referee proposes to use
\cite[Thm.~5.1]{Chapman};
he thus made Theorem~\ref{reembedding} stronger:
now $h$ is a self-homeomorphism of $\mathbb R^n$ (not only an embedding of $X$, as I stated in the preprint).


\begin{thebibliography}{70}


\bibitem{Alexander-sphere}
J.W. Alexander, An example of a simply connected
surface bounding a region which is not
simply connected,
Proc. Nat. Acad. Sci. USA 10 (1) (1924) 8--10.

\bibitem{Alexander}
J.W. Alexander,
Remarks on a point set constructed by Antoine,
Proc. Nat. Acad. Sci. USA 10 (1) (1924) 10--12.

\bibitem{Antoine}
L. Antoine,
Sur la possibilite d'etendre l'hom\'eomorphie de deux
figures a leur voisinages,
C.R. Acad. Sci. Paris 171 (1920) 661--663.

\bibitem{Antoine-173}
L. Antoine,
Sur les ensembles parfaits partout discontinus, C.R. Acad. Sci. Paris 173 (1921) 284--285.

\bibitem{Antoine-diss}
L. Antoine,
Sur l'hom\'eomorphie de deux
figures
et de  leurs voisinages.
Th\`eses de l'entre-deux-guerres,
vol. 28, 1921.
http://eudml.org/doc/192716


\bibitem{Antoine-FM}
L. Antoine,
Sur les voisinages de deux figures hom\'eomorphes,
Fund. Math. 5 (1924) 265--287.

\bibitem{Apodaca}
E.R. Apodaca,
On the simultaneous embedding of uncountably many distinct wild arcs with one wild endpoint in $E^3$, a geometric approach,
Fund. Math. 113 (1981) 175--186.

\bibitem{BEM}
H. Becker, F. van Engelen, J. van Mill,
Disjoint embeddings of compacta,
Mathematika 41 (2) (1994) 221--232.


\bibitem{Bing52}
R.H. Bing,
A homeomorphism between the $3$-sphere
and the sum of two solid horned spheres,
Ann. of Math. 56 (2) (1952) 354--362.


\bibitem{Bing-abstr}
R.H. Bing,
$E^3$ does not contain
uncountably many mutually exclusive wild
surfaces, 
Bull. Amer. Math. Soc.
63 (1957) 404.
Abstract \# 801t.


\bibitem{Bing59}
R.H. Bing,
Conditions under which a surface in $E^3$ is tame,
Fund. Math. 47 (1959) 105--139.

\bibitem{Bing-tame}
R.H. Bing,
Tame Cantor sets in $E^3$,
Pacif. J. Math. 11 (2) (1961) 435--446.

\bibitem{Bing-TAMS61}
R.H.~Bing,
A surface is tame if its complement is 1-ULC,
Trans. Amer. Math. Soc. 101 (1961), 2, 294--305.

\bibitem{Blankinship}
W.A. Blankinship,
Generalization of a construction of Antoine,
Ann. Math. (2) 53 (1951) 276--297.

\bibitem{Bogachev}
V.I. Bogachev.
Measure Theory.
Volume II.
Springer, 2007.

\bibitem{BM}
B.L. Brechner, J.C. Mayer, 
Antoine's Necklace or How to Keep a Necklace From Falling
Apart,
Coll. Math. J. 19 (4) (1988) 306--320.

\bibitem{Bryant}
J.L. Bryant,
Concerning uncountable families of $n$-cells in $E^n$,
Michigan Math. J. 15 (1968) 477--479.

\bibitem{BC}
C.E. Burgess, J.W. Cannon,
Embeddings of surfaces in $E^3$,
Rocky Mountain J. of Math.
1 (2) (1971) 259--344.

\bibitem{Burgess-75}
C.E.~Burgess,
Embeddings of surfaces in Euclidean three-space, 
Bull. Amer. Math. Soc.
81 (1975), 5, 795--818.

\bibitem{Chapman}
T.A. Chapman. Lectures on Hilbert Cube Manifolds.
Expository lectures from the CBMS Regional Conference held at Guilford College, October 11-15, 1975. 
Conference Board of the Mathematical Sciences Regional Conference Series in Mathematics. No. 28. Providence, R.I.: American Mathematical Society. 1976.

\bibitem{Cernavsky}
A.V. Chernavskii,
On the work of L.V. Keldysh and her seminar,
Russian Math. Surveys, 60 (4) (2005) 589--614;
transl. from: Usp. Mat. Nauk 60 (4) (2005)  11--36.

\bibitem{Daverman-absense}
R.J. Daverman,
On the absence of tame disks in some wild cells,
Geom. Topol., Proc. Conf. Park City 1974, Lect. Notes Math. 438, 1975. 142--155.

\bibitem{Daverman}
R.J.~Daverman, 
Embeddings of $(n-1)$-spheres in
Euclidean $n$-space,
Bull. Amer. Math. Soc.
84 (3) (1978) 377--405.

\bibitem{D-dec}
R.J. Daverman.
Decompositions of Manifolds. AMS Chelsea Publishing.
American Mathematical Society. 
Providence, Rhode Island, 2007.

\bibitem{DV}
R.J. Daverman, D.A. Venema.
Embeddings in Manifolds.
Graduate Studies in Mathematics 106. Providence, RI: American Mathematical Society.
2009.

\bibitem{D1}
A. Denjoy,
Sur les ensembles parfaits discontinus,
C.R. Acad. Sci. Paris 149 (1909) 1048--1050.

\bibitem{D2}
A. Denjoy,
Continu et discontinu,
C.R. Acad. Sci. Paris 151 (1910) 138--140.

\bibitem{vanDouwen}
E.K. van Douwen,
Uncountably many pairwise disjoint copies
of one metrizable compactum in another,
Topol. Appl. 51 (1993) 87--91.

\bibitem{Ea}
W.T. Eaton,
A Generalization of the Dog Bone Space to $E^n$,
Proc. Amer. Math. Soc. 39 (2) (1973) 379--387.

\bibitem{Engelking}
R. Engelking. Dimension Theory. 
PWN Polish Scientific Publishers/North-Holland Publishing Company, Warszawa/Amsterdam, 1978.

\bibitem{Frolkina-thesis}
O. Frolkina,
Uncountable families of pairwise disjoint
inequivalent wild $k$-disks in $\mathbb R^n$, 
Abstracts of the Conference
``Problems of Modern Topology and its
Applications'',
11--12 May 2017, Tashkent, Uzbekistan.
82--83.

\bibitem{GRZ2005}
D. Garity, D. Repov\v{s}, M. \v{Z}eljko,
Uncountably many inequivalent
Lipschitz homogeneous Cantor sets in 
$ \mathbb R^3$,
Pacif. J. Math. 222 (2) (2005) 287--299.

\bibitem{GO}
D.G. DeGryse, R.P. Osborne,
A wild Cantor set in $E^n$ with simply connected complement,
Fund. Math. 86 (1974) 9--27.

\bibitem{Ivanov}
A.A. Ivanov,
Isotopy of compacta in Euclidean spaces,
Soviet Math. Dokl. (=Dokl. Akad. Nauk SSSR)
71 (6) (1950) 1021--1022 (in Russian).
See Zbl 0037.26305.

\bibitem{Ivanov-diss}
A.A. Ivanov,
Isotopy of compacta in Euclidean spaces,
Doctoral dissertation 
(Candidate of Sciences),
St. Petersburg Department of
V.A. Steklov Institute of Mathematics of
the USSR Academy of Sciences,
Moscow-Leningrad, 1950
(in Russian).

\bibitem{Keldysh62}
L.V. Keldysh,
Embedding of zero-dimensional
compact sets in $E^n$,
Dokl. Akad. Nauk SSSR 147 (1962) 772--775;
English transl. in: Soviet Math. Dokl.
3 (1962).

\bibitem{Keldysh}
L.V. Keldysh,
Topological imbeddings in Euclidean space (English. Russian original),
Proc. Steklov Inst. Math. 81 (1966), 203 p.; translation from: Tr. Mat. Inst. Steklov. 81  (1966), 184 p.

\bibitem{Martin}
J. Martin,
A note on the uncountably many disks,
Pacific J. Math. 13 (1963) 1331--1333.

\bibitem{McMillan-Taming}
D.R. McMillan, Jr., Taming Cantor sets in $E^n$,
Bull. Amer. Math. Soc. 70 (1964) 706--708.

\bibitem{Me}
S. Melikhov,
Review of ``Embeddings in manifolds''
by R.J.~Daverman and G.A.~Venema,
Mathematical Reviews,
MR2561389 (2011g:57025).

\bibitem{Moise}
E.E. Moise.
Geometric Topology in Dimensions $2$ and $3$.
Springer-Verlag, 1977.

\bibitem{Osborne1}
R.P. Osborne,
Embedding Cantor sets in a manifold,
Part I: Tame Cantor sets in $E^n$,
Michigan Math. J. 13 (1966) 57--63.

\bibitem{Osborne2}
R.P. Osborne,
Embedding Cantor sets in a manifold,
Part II: An extension theorem for homeomorphisms on Cantor sets,
Fund. Math. 65 (2) (1969) 147--151.

\bibitem{Osborne3}
R.P. Osborne,
Embedding Cantor sets in a manifold,
Part III: Approximating spheres,
Fund. Math. 90 (3) (1976) 253--259.


\bibitem{Rushing}
T.B. Rushing. Topological embeddings.
Pure and Applied Mathematics, Vol.~52,
New York and London, 
Academic Press, 1973.

\bibitem{Sher-note}
R.B. Sher,
A note on the example of Stallings,
Proc. Amer. Math. Soc., 19 (3) (1968) 619--620.

\bibitem{Sher68}
R.B. Sher, Concerning wild Cantor sets in $E^3$,
Proc. Amer. Math. Soc. 19 (5) (1968) 1195--1200.

\bibitem{Sher69}
R.B. Sher, 
Families of arcs in $E^3$,
Trans. Amer. Math. Soc. 143 (1969) 109--116.

\bibitem{Shtanko69}
M.A. Shtan'ko,
The imbedding of compacta in Euclidean space, Soviet Math. Dokl. 10 (1969) 758--761.
Transl. from: Dokl. Akad. Nauk SSSR 186 (1969) 1269-1272.

\bibitem{Shtanko}
M.A. Shtan'ko, The embedding of compacta in Euclidean space, Math. USSR-Sb. 12 (2) (1970) 234--254.

\bibitem{Stallings}
J.R. Stallings, 
Uncountably many wild disks,
Annals of Math. 71 (1) (1960) 185--186.

\bibitem{TW}
E.D. Tymchatyn, R.B. Walker,
Taming the Cantor fence,
Topol. Appl. 83 (1998) 45--52.

\bibitem{Urysohn}
P. Urysohn,
M\'emoire sur les multiplicit\'es
Cantoriennes,
Fund. Math. 7 (1) (1925) 30--137.

\bibitem{Urysohn1}
P.S. Urysohn,
Works in topology
and other fields of mathematics. Vol. 1, GITTL, Moscow, 1951. (In Russian.)

\bibitem{Wright1983}
D.G. Wright,
Ambiently universal sets in $E^n$,
Trans. Amer. Math. Soc. 277 (2) (1983) 655--664.

\bibitem{Wright1986}
D.G. Wright.
Rigid sets in manifolds. In:
Geom. and Alg. Topology.
Banach Center Publ.
18 (1986) 161--164.

\bibitem{Wright-rigid}
D.G. Wright,
Rigid sets in $E^n$,
Pacif. J. Math. 121 (1) (1986) 245--256.

\bibitem{Z}
M. \v{Z}eljko,
Minimal number of tori in geometric
self-similar Antoine Cantor sets, 
JP J. Geom. Topol. 5 (2) (2005) 109--113.

\bibitem{ZS}
G.Ya. Zyev, A.B. Sosinsky,
On spheres containing wild zero-dimensional sets,
Moscow Univ. Math. Bull. 3 (1967) 5--9 (in Russian).

\end{thebibliography}
\end{document}